\documentclass[12pt]{amsart}

\usepackage{amscd}
\usepackage{amsthm}
\usepackage{amsfonts}
\usepackage{amssymb}
\usepackage{amsmath}
\usepackage{amsxtra}
\usepackage{eucal}

\newcommand{\nc}{\newcommand}
\nc{\C}{{\mathcal C}}
\nc{\CB}{{\mathcal B}}
\nc{\CO}{{\mathcal O}}
\nc{\CK}{{\mathcal K}}
\nc{\Spec}{{\operatorname{Spec}}}
\nc{\CG}{{\mathcal G}}
\nc{\CA}{{\mathcal A}}
\nc{\CU}{{\mathcal U}}
\nc{\supp}{\operatorname{supp}}
\nc{\LL}{{\mathbf L}}

\nc{\D}{{\mathbb D}}
\nc{\St}{\operatorname{St}^{\bullet}}
\nc{\CM}{{\mathcal M}}
\nc{\A}{{\mathcal A}}

\nc{\1}{{\mathbf 1}}
\nc{\CC}{{\mathbb C}}
\nc{\CR}{{\mathcal R}}
\renewcommand{\k}{{\mathbb  k}}
\nc{\Q}{{\mathbb Q}}
\nc{\M}{{\mathcal M}}
\nc{\CL}{{\mathcal L}}
\nc{\U}{{\mathbf U}}
\nc{\B}{{\mathbf B}}
\nc{\Z}{{\mathbb Z}}
\nc{\CRhom}{\operatorname{RHom}\bul}
\nc{\Ad}{\operatorname{Ad}}
\nc{\CRes}{\operatorname{Res}}
\nc{\gr}{\operatorname{gr}}
\nc{\tr}{\operatorname{tr}}
\nc{\End}{\operatorname{End}}

\nc{\g}{{\mathfrak g}}
\nc{\hatg}{\hat{\frak g}}

\renewcommand{\b}{{\mathfrak b}}
\nc{\sem}{{$\S_{\g}^{\g_{>0}}}}
\nc{\gl}{{\frak g\frak l}}
\nc{\n}{{\frak n}}
\nc{\si}{{\frac\infty 2}}
\nc{\p}{{\frak p}}
\nc{\h}{{\frak h}}

\nc{\Ind}{\operatorname{Ind}}
\nc{\ch}{\operatorname{ch}}
\nc{\Coind}{\operatorname{Coind}}
\nc{\opp}{{\operatorname{opp}}}
\nc{\Ker}{\operatorname{Ker}}
\nc{\im}{\operatorname{Im}}
\nc{\Coker}{\operatorname{Coker}}
\nc{\dirlim}{\underset{\rightarrow}{\operatorname{lim}}}
\nc{\invlim}{\underset{\leftarrow}{\operatorname{lim}}}
\nc{\Sem}{{{\mathbf S}_{\g}^{\g_{>0}}}}
\nc{\CN}{{\mathcal N}}
\nc{\Ext}{\operatorname{Ext}^{\bullet}}
\nc{\ext}{\operatorname{Ext}}
\nc{\tilW}{\til{W}}
\nc{\Mat}{\mathcal{M}at}
\nc{\CV}{\mathcal{V}}
\nc{\lth}{\ell t}
\nc{\BB}{\mathcal{B}}
\nc{\Tor}{\operatorname{Tor}_{\bullet}}
\nc{\tor}{\operatorname{Tor}}
\nc{\Tors}{\operatorname{Tor}_{\frac \infty 2+\bullet}}
\nc{\Exts}{\operatorname{Ext}^{\frac \infty 2+\bullet}}
\nc{\Hom}{\operatorname{Hom}^{\bullet}}
\nc{\ad}{\operatorname{ad}}
\renewcommand{\hom}{\operatorname{Hom}}
\nc{\Tate}{_{\mathsf{Tate}}}
\renewcommand{\mod}{\operatorname{mod}}
\nc{\Mod}{\operatorname{Mod}}

\nc{\Barb}{\operatorname{Bar}^{\bullet}}

\nc{\upX}{X^{\uparrow}}
\nc{\upcD}{{\mathcal D}^{\uparrow}}
\nc{\upD}{D^{\uparrow}}
\nc{\dX}{X^{\downarrow}}
\nc{\dcD}{{\mathcal D}^{\downarrow}}
\nc{\dD}{D^{\downarrow}}
\nc{\upC}{{\mathcal C}^{\uparrow}}
\nc{\dC}{{\mathcal C}^{\downarrow}}
\nc{\underA}{\underline{A}}
\nc{\underC}{\underline{\CC}}
\nc{\underB}{\underline{B}}
\nc{\underk}{\underline{\k}}
\nc{\Db}{D^{\bullet}}
\nc{\ten}{{\otimes}}
\nc{\tenb}{{\boxtimes}}
\nc{\tenf}{\overset{\operatorname{!}}\ten}
\nc{\tenl}{\overset{\operatorname{L}}\ten}
\nc{\map}{\longrightarrow}
\nc{\eps}{\varepsilon}
\nc{\bs}{\bigskip\\}
\nc{\ms}{\smallskip\\}

\nc{\tilbar}{\widetilde{\operatorname{Bar}}}
\nc{\tilBarb}{\widetilde{\operatorname{Bar}}^{\bullet}}
\nc{\overr}{\overline{R}}
\nc{\overI}{\overline{I}}
\nc{\overX}{\overline{X}}
\nc{\overh}{\overline{h}}
\nc{\overY}{\overline{Y}}
\nc{\overW}{\overline{W}}
\nc{\linbar}{\overline{\operatorname{Bar}}}

\nc{\til}{\widetilde}
\nc{\oppA}{A^{\sharp}}
\nc{\Lemma}{{\bf Lemma:\ }}
\nc{\Theorem}{{\bf Theorem:}\ }
\nc{\Cor}{{\bf Corollary:}\ }
\nc{\Def}{{\bf Definition:}\ }
\nc{\Prop}{{\bf Proposition:}\ }
\nc{\Con}{{\bf Conjecture:}\ }
\nc{\Rem}{{\bf Remark:}\ }
\nc{\dok}{{\bf Proof.}\ }

\nc{\SInd}{\operatorname{S-Ind}}
\nc{\SCoind}{\operatorname{S-Coind}}
\nc{\bul}{^{\bullet}}
\nc{\stand}{C^{\frac\infty2+\bullet}}
\nc{\ssn}{\subsection{}}
\nc{\sssn}{\subsubsection{}}

\nc{\hgt}{\operatorname{ht}}
\nc{\sqbinom}{\fracwithdelims[][0pt]}

\address{Independent University of Moscow, Pervomajskaya st. 16-18,
Moscow 105037, Russia}
\email{hippie@mccme.ru}
\author{Sergey Arkhipov}
\title{Semiinfinite cohomology of  Lie-* algebras}
\date{}

\begin{document} 
\maketitle 
\section{Introduction.}  
This paper is a
natural extension of the previous note \cite{Ar2}.  Semiinfinite
cohomology of Tate Lie algebra was defined in that note in terms of
some duality resembling Koszul duality.  The language of differential
graded Lie algebroids was the main technical tool of the note.

The present note is devoted to globalization of the main construction
from \cite{Ar2} in the following sense.  The setup in \cite{Ar2}
included a suitably chosen module over a Tate Lie algebra $\g$ with a
fixed Lie subalgebra $\b$ being a {c-lattice} in $\g$.

The rough global analogue of this picture is as follows.  Consider a
compact curve $X$ over a field of characteristic zero.  Denote by
$\operatorname{mod-}D_X$ the category of (right) D-modules on $X$.  We
fix a Lie algebra $\CG$ in the category $\operatorname{mod-}D_X$.
This data can be viewed as a family of Lie algebras $\g_x$ along the
curve $X$.  Another part of the data includes a Lie subalgebra
$\CB\subset\CG$.  So the problem is to define semiinfinite cohomology
complex of such pair.

In fact we need some additional constraints on the pair
$\CB\subset\CG$.  So the formal picture starts from a different notion
of a {Lie-* algebra} $\CL$ on $X$ (see the precise definition in
Section 2).  Roughly spaking a $D_X$-locally free Lie-* algebra is a
D-module incarnation of a Lie algebra in the category of vector
bundles on $X$ with the bracket given by a differential operator.  We
define two types of modules over a Lie-* algebra (see \ref{modules}).
The first one called a {\em Lie-* module} is just a D-module
incarnation of the module over the Lie algebra in the category of
vector bundles, like above, with the action given by a differential
operator.  Still we will be more interested in the second type of
modules over a Lie-* algebra called {\em chiral modules} (see
\ref{modules} for the definition).

So starting from a Lie-* algebra $\CL$ and a chiral module $\CM$ we
perform the main construction more or less parallel to the one from
\cite{Ar2}.  Namely we define the Lie algebra $\CG=\CG(\CL)$ in the
category $\operatorname{mod-}D_X$ with the Lie subalgebra
$\CB=\CB(\CL)\subset\CG(\CL)$.  We show that a $\CL$-chiral module
$\CM$ becomes a $\CG(\CL)$-module.

Next, imitating the construction of \cite{Ar2} Section 4, we define a
DG Lie algebroid $\CA\bul(\CL)$ in the category of $D_X$-modules over
a DG $\tenf$-algebra $\CR\bul=\CR\bul(\CL)$.  Koszul duality type
construction provides a {\em left} DG-module $C\bul(\CL,\CM)$ over
$\CA\bul(\CL)$.

To go further one needs to pass to a central extension of $\CL$ called
the Tate central extension and denoted by $\CL_{\Tate}$.  Let $\CM$ be
a chiral module over $\CL_{\Tate}$.  It turns out that the complex of
D-modules $C\bul(\CL_{\Tate},\CM)$ by some antipode construction
becomes a {\em right} module over $\CA\bul(\CL)$.

Finally we consider the homological Chevalley complex of the DG Lie
algebroid $\CA\bul(\CL)$ in the category of $D_X$-modules with
coefficients in $C\bul(\CL_{\Tate},\CM)$ (see \ref{homol} for the
definition of the homological Chevalley complex of a DG Lie
algebroid).  We call the obtained complex of D-modules the global
standard semiinfinite complex of the Lie-* algebra $\CL$ with
coefficients in the chiral module $\CM$.

Let us say a few words about the structure of the paper.  In Section 2
we collect necessary definitions and simple facts about Lie-*
algebras, Lie-* modules, chiral modules etc.  Section 3 is devoted to
the construction of the Tate central extension of a Lie-* algebra.
Section 4 contains all the necessary constructions concerning DG Lie
algebroids in the category of $D_X$-modules.  In particular we present
the definition of the homological Chevalley complex of of a DG Lie
algebroid with coefficients in a {\em right} DG-module.  Section 5 is
the heart of the paper.  We present the constructions of the Lie
algebras $\CB(\CL)$ and $\CG(\CL)$ in the category of $D_X$-modules.
Then we define the DG Lie algebroid $\CA\bul(\CL)$ over the DG
$\tenf$-algebra $\CR\bul(\CL)$.  Finally after overcoming the problem
of necessity to pass to the Tate central extension of $\CL$ we present
the standard semiinfinite complex $C^{\si}(\CL_{\Tate},\CM)$ for a
$\CL_{\Tate}$-chiral module $\CM$.

Note that the paper \cite{BD} contains a construction of the global
BRST complex for calculating the semiinfinite cohomology of a chiral
module over a Lie-* algebra.  Somehow the present paper grew out of an
attempt to understand that construction avoiding the notions of chiral
algebras, chiral enveloping algebras etc.The technique used in the
definition in \cite{BD} is quite different from ours and it is not
checked that the two constructions give the same answer.

{\bf Acknowledgements.}  The author is happy to thank Sasha Beilinson
and Dennis Gaitsgory who explained him the chiral algebra basics.  The
author also would like to express his deep gratitude to IAS,
Princeton, USA, and IHES, Bures-sur-Yvette, France, where parts of the
work on the paper were done, for hospitality and extremely stimulating
working conditions.

\section{Lie-* algebras.}  
In this section we recall briefly basic
notation and constructions concerning Lie-* algebras.  In our
exposition we follow \cite{Denis}.  We will be working over a fixed
smooth curve $X$.  Denote the diagonal embedding $X\hookrightarrow
X\times X$ by $\Delta$.  The embedding of the complementary open set
$X\times X\setminus X_\Delta\hookrightarrow X\times X$ is denoted by
$j$.

Let $D_X\operatorname{-}\mod$ (resp.  $\mod\operatorname{-}D_X$) be
the category of left (resp.  right) modules $\CM$ over the sheaf of
algebras of differential operators on $X$ such that $\CM$ is
quasicoherent over $\CO_X$.  It is known that the category
$D_X\operatorname{-}\mod$ is naturally a symmetric tensor category
with the tensor product given by $\CM\ten_{\CO_X}\CN$ for $\CM, \CN\in
D_X\operatorname{-}\mod$.  Let $\Omega=\Omega^1_X$.  Then the category
$\mod\operatorname{-}D_X$ becomes a symmetric tensor category with the
tensor product given by 
$$
\CM\tenf\CN=(\CM\ten_{\CO_X}\Omega^{-1})\ten_{\CO_X}(\CN\ten_{\CO_X}
\Omega^{-1})\ten_{\CO_X}\Omega.  
$$ 
\subsection{Definition.}  Recall
that a {\em Lie-* algebra} on $X$ is a right D-module $\CL$ with the
map 
$$ 
\{\cdot\tenb\cdot\}:\ \CL\boxtimes\CL\map\Delta_!(\CL) 
$$ 
which
is antisymmetric and satisfies the Jacobi identity in the following
sense.  If $a\boxtimes b\boxtimes c\cdot f(x,y,z)$ is a section of the
$\CL\boxtimes \CL\boxtimes\CL$ on $X\times X\times X$, then the
element 
$$
\{\{f(x,y,z)\cdot a\boxtimes b\}\boxtimes c\}+
\sigma_{1,2,3}(\{\{f(z,x,y)\cdot b\boxtimes c\}\boxtimes a\})+
\sigma_{1,2,3}(\{\{f(y,z,x)\cdot c\boxtimes a\}\boxtimes b\})
$$ 
of
$\Delta_{x=y=z}{}_{!}(\CL)$ vanishes.  Here $\sigma_{1,2,3}$ denotes
the lift of the cyclic automorphism of $X\times X\times X$:
$(x,y,z)\mapsto (y,z,x)$ to the D-module $\Delta_{x=y=z}{}_{!}(\CL)$.

\sssn 
Note that if $\CL$ is a Lie-* algebra, it follows from the
definition that $DR^0(\CL)$ is a sheaf of ordinary Lie algebras;
moreover it acts on $\CL$ by endomorphisms of the D-module structure
that are derivations of the Lie-* structure.

In particular, for an affine subset $U\subset X$, $DR^0(U,\CL)$ is a
Lie algebra.  Thus for any point $x\in X$ the topological spaces
$DR^0(\Spec(\hat\CO_x),\CL)$ and $DR^0(\Spec(\hat\CK_x),\CL)$ carry
the natural structures of topological Lie algebras.

\sssn 
Our next step consists of reformulation of the Lie$^*$ algebra
definition in terms of Lie coalgebras in the standard tensor structure
$\overset{!}{\ten}$ on right D-modules.
\vskip 1mm
\noindent
\Lemma \label{hom}
Let $\CM_1$, $\CM_2$ be two D-modules on $X$ with $\CM_1$ being
locally free and finitely generated. Then:
\begin{itemize}
\item[(i)]
For a third D-module $\CM$ on $X$, there is a canonical
isomorphism:
$$\hom_{D_X^2}(\CM\boxtimes \CM_1,\Delta_{!}(\CM_2))\til\map
\hom_{D_X}(\CM,\CM_1^*\overset{!}\otimes\CM_2).$$
\item[(ii)]
The canonical map (from point (i))
$\left(\CM_1^*\overset{!}\otimes\CM_2\right)
\boxtimes\CM_1\map\Delta_{!}(\CM_2)$ induces
an isomorphism
$DR^0(\CM_1^*\overset{!}\otimes\CM_2)\map \hom_D(\CM_1,\CM_2)$.
\qed
\end{itemize}
Below we always suppose that any Lie-* algebra $\CL$ we work with is
locally free and finitely generated as a D-module on $X$.
\vskip 1mm
\noindent
\Cor
For a Lie-* algebra $\CL$ we have canonical maps
$$\operatorname{co-ad}:\CL\boxtimes\CL^*\map\Delta_{!}(\CL^*) \text{ and }
\operatorname{co-br}:\CL^*\to \CL^*\overset{!}\otimes\CL^*.\qed
$$
We shall call the two maps of the above corollary ``the co-adjoint
action'' and ``the co-bracket'', respectively.

\sssn
\Rem  \label{coalgebra}
In particular note that to specify a structure
of a Lie-* algebra on a D-module $\CL$ is the same as to
provide a structure of a {\em Lie coalgebra (in the usual tensor
structure on the category of right D-modules)} on the D-module
$\CL^*:=\underline{\hom}_{D_X}(\CL,D_X\ten\Omega_X)$. 
The co-bracket here is the one
obtained in the previous Corollary.

\subsection{Modules over a Lie-* algebra.}
There are two different ways to define a notion of a module over a Lie-*
algebra.

\sssn                                       \label{modules}
{\em A Lie-* module} over a Lie-* algebra $\CL$ is a (right) D-module
$\CM$ on $X$ with a map $$\rho:\CL\boxtimes\CM\to\Delta_{!}(\CM)$$ such that
for a section $a\boxtimes b\boxtimes m$ of $\CL\boxtimes\CL\boxtimes\CM$ the
two sections $$\rho(\{a\boxtimes b\},m) \text{ and
}\rho(a,\rho(b,m))-\sigma_{1,2}(\rho(b,\rho(a,m)))$$ of
$\Delta_{1,2,3}{}_{!}(\CL)$ coincide.

\smallskip
 
Like in  the Lie-* algebra case a structure of the Lie-* module on 
$\CM$ provides a structure of the sheaf of modules over the sheaf of Lie 
algebras
$DR^0(\CL)$ on $DR^0(\CM)$. Moreover the sheaf of Lie algebras $DR^0(\CL)$
acts on $\CM$ and contrarywize such an action recovers a Lie-* module 
structure
if it is given by differential operators (see \cite{BD}, 2.5.4).
\vskip 1mm
\noindent
\Rem
Similarly to \ref{coalgebra} one can easily check that to specify a
structure of a Lie-* module over the  Lie-* algebra $\CL$ on a
D-module $\CM$ is the same as to make $\CM$ a comodule over the Lie
coalgebra $\CL^*$, i.e. to provide a D-module map
$\operatorname{co-ac}:\ \CM\map\CL^*\tenf\CM$
satisfying certain coassociativity constraints.

\sssn
{\em A chiral module} over a Lie-* algebra is again a (right)
D-module $\CM$ on $X$, but
with an operation $$\rho:j_*j^*(\CL\boxtimes\CM)\to\Delta_{!}(\CM)$$ such
that every section $f(x,y,z)\cdot a\boxtimes b\boxtimes
m\in\Gamma(X\times X\times X\setminus (\Delta_{x=z}\cup \Delta_{y=z}),
\CL\boxtimes\CL\boxtimes\CM)$ satisfies an identity as follows:
$$\rho(\{f(x,y,z)\cdot a\boxtimes b\},m)= \rho(a,\rho(f(x,y,z)\cdot
b,m))-\sigma_{1,2}(\rho(b,\rho(f(y,x,z)\cdot a,m)),$$ as sections of
$\Delta_{x=y=z}{}_{!}(\CM)$.

\subsection{Naive de Rham complex for a Lie-* module over a Lie-*
algebra.} Let $\CM$ be a Lie-* module over a (locally D-free finitely
generated) Lie-* algebra $\CL$. Consider a complex $C\bul(\CL,\CM)$ of
D-modules on X as follows:
$C^k(\CL,\CM)=\Lambda^k(\CL^*)\tenf \CM$ as a D-module. The differential
is given by
$$d(\omega_1\wedge\ldots\wedge\omega_k\ten m)
=\sum \omega_1\wedge\ldots\wedge\operatorname{co-br}(\omega_i)
\wedge\ldots\wedge\omega_k\ten m
+\omega_1\wedge\ldots\wedge\omega_k\wedge\operatorname{co-ac}(m).
$$
Here $\omega_i$ are the sections of $\CL^*$ and $m$ is the section of $\CM$.
\vskip 1mm
\noindent
\Lemma
The complex of D-modules $C\bul(\CL,\CM)$ is well-defined, i.e. $d^2=0$.
\qed
\vskip 1mm
\noindent
\Rem
We call the constructed complex the {\em naive de Rham complex} for
the Lie-* module over the Lie-* algebra.  However below we construct a
more sophisticated complex for a {\em chiral} module over a Lie-*
algebra. That complex is also of de Rham origin. It will be realized
as  de Rham complex of a certain DG Lie-* module over a certain {\em DG
Lie-* algebroid}.

\subsection{De Rham DG-algebra of a Lie-* algebra.}
Consider the complex $C\bul(\CL,\CM)$ for $\CM$ equal to
the trivial Lie-* module $\Omega$ over $\CL$.
\vskip 1mm
\noindent
\Lemma
The wedge product makes the complex $C\bul(\CL,\Omega)$ into a
supercommutative DG-algebra in the tensor category of right D-modules.
\qed

We denote this DG-algebra by $\CR\bul=\CR\bul(\CL)$.

\section{Tate extension of a Lie-* algebra.}
\subsection{Matrix Lie-* algebra.}
Let $\CV$ be a locally free finitely generated D-module on $X$, and
let $\CV^*:=\underline{\hom}_{D_X}(\CV,D_X\ten\Omega_X)$ be its
(Verdier) dual D-module.  Consider the $\tenf$-tensor product of $\CV$
and $\CV^*$.

\sssn
\Rem
By Lemma \ref{hom}(i) we have a Lie-* pairing 
$$\langle\cdot\tenb\cdot\rangle:\ \CV\tenb\CV^*\map\Delta_!(\Omega).$$
\Lemma
The D-module $\CV\tenf\CV^*=:\Mat$ carries a natural structure of an 
associative-* algebra on $X$.

\dok
The associative product map is defined as follows:
\begin{multline*}
\operatorname{as}:\ \left(\CV_{(1)}\tenf\CV_{(1)}^*\right)
\tenb\left(\CV_{(2)}\tenf\CV_{(2)}^*\right)\\
\til\map
\left((\CV_{(1)}\tenb\Omega)\tenf(\Omega\tenb\CV_{(2)}^*)\right)
\tenf\left((\CV_{(2)}\tenb\Omega)\tenf(\Omega\tenb\CV_{(1)}^*)\right)\\
\overset{\langle\cdot\tenb\cdot\rangle\tenf\operatorname{id}}{\map}
\Delta_!(\Omega)\tenf\left((\CV_{(2)}\tenb\Omega)\tenf(\Omega\tenb\CV_{(1)}^*)
\right)\til\map
\Delta_!(\CV\tenf\CV^*).
\end{multline*}
Here the indexes $(1)$ and $(2)$ denote the factors in the product,
and  the *-pairing is taken between the ones in the first pair of brackets.
It is left to the reader to check the associativity of the product.\qed

Now we obtain the Lie-* bracket on $\Mat$ from the associative-*
product in the usual way:  $$a\tenb b\mapsto \operatorname{as}(a\tenb
b)-\operatorname{as}(b\tenb a).$$ We call $\Mat$ the {\em matrix Lie-*
algebra} of the D-module $\CV$.

Note that the associative-* algebra $\Mat$ acts both on $\CV$ and on
$\CV^*$ in a canonical way and the pairing
$\langle\cdot\tenb\cdot\rangle$ is $\Mat$-invariant.

\subsection{Tate extension of the Lie-* algebra $\Mat$.}
The material of this subsection is almost word to word copied from
\cite{BD}, 2.6. We include it in our paper just for the sake of 
completeness.

\sssn
The Tate extension is a
canonical central extension of
Lie-* algebras
$$
0\map\Omega_X\map \Mat^\flat
\overset{\pi^\flat}\map
\Mat\map 0. 
$$

To define $\Mat^\flat$ as a
D-module consider the exact sequence of
D-modules on
$X\times X$
$$
\CV\boxtimes \CV^*\overset{\eps}\map j_*j^*(\CV\boxtimes \CV')
\overset{\pi}\map\Delta_!(\CV\tenf \CV^*)\map 0.
$$
Here as before $\Delta:\ X\hookrightarrow X\times X$
is the diagonal embedding, 
$j:\ U:=X\times X\smallsetminus
\Delta(X)\hookrightarrow X\times X$ is the 
complementary
embedding, and $\pi$ is the canonical
arrow. 

Namely one has
$\CV\tenf\CV^*=H^1\Delta^!(\CV\boxtimes \CV^*)=
\Coker
\,\varepsilon$; explicitly, 
$\pi$ sends
$(t_2-t_1)^{-1}v\boxtimes v'\in
j_*j^*\CV\boxtimes \CV^*$ to $v\otimes
v'(dt)^{-1}\in\Delta_\cdot(\CV\tenf\CV')\subset\Delta_!(\CV\tenf\CV^*)$. 
Note that $\langle\cdot\tenb\cdot\rangle$ vanishes on $\Ker\,\eps$ and,
pushing forward the above exact sequence by
$\langle\cdot\tenb\cdot\rangle$, we get an extension of
$\Delta_!(\CV\tenf \CV^*)$ by $\Delta_!(\Omega)$.

This extension is supported on the diagonal. Applying
$\Delta^!$ we get the Tate extension
$\Mat^\flat$. We denote the canonical morphism
$j_*j^*(\CV\boxtimes \CV^*)\map
\Delta_! \Mat^\flat$ as $\mu_{Mat}$.

\sssn
The above construction is natural with respect
to Lie-* algebras actions. Namely, assume
that a Lie-* algebra
$\CL$ acts on $\CV$,
$\CV^*$ so that $\langle\cdot\tenb\cdot\rangle$ is
$\CL$-invariant. Then the $\CL$-action on
$\Mat =\CV\tenf \CV^*$ lifts canonically to
an $\CL$-action on $\Mat^\flat$. To see this,
consider the $D_X^{\otimes 2}$-modules
$\Delta^\cdot(\CV\boxtimes \CV^*)$, $\Delta^\cdot
j_*j^*(\CV\boxtimes \CV^*)$. The Lie algebra $h(\CL)$
acts on them in the obvious manner. The
morphism $\langle\cdot\tenb\cdot\rangle:\ \Delta^\cdot(\CV\boxtimes
\CV^*)\map \Delta^\cdot\Delta_!(\Omega)$ is
$DR^0(\CL )$-invariant. Therefore $DR^0(\CL)$ acts on
$\Mat^\flat$. This action is uniquely
 determined by property that
$\mu_{\Mat}$ is a morphism of
$DR^0(\Mat)$-modules.
Evidently the action of $DR^0(\CL)$
is given by differential operators
so we have the desired
$\CL$-action on $\Mat^\flat$.

In particular, the canonical $\Mat$-actions
on $\CV$, $\CV^*$ define a $\Mat$-action on
$\Mat^\flat$ that lifts the adjoint action on
$\Mat$. Composing this action with
$\pi^\flat$ we get an operation  $$\{\cdot\tenb\cdot\}:\
\Mat^\flat\tenb\Mat^\flat\map\Delta_!(\Mat^\flat ).$$

\sssn 
\Lemma
This is a Lie-* bracket.
\qed

\subsection{Tate extension of a Lie-* algebra.}
Now let $\CL$ be an arbitrary $D_X$-locally free 
finitely generated Lie-* algebra.
Recall that by Lemma~\ref{hom} we have a canonical co-action map of 
the $D_X$-modules
$$\operatorname{co-ac}:\CL\map \CL\overset{!}\otimes\CL^*.$$
We interprete as a (D-module) map $\operatorname{can}:\ \CL\map\Mat$.
\vskip 1mm
\noindent
\Lemma
The map $\operatorname{can}$ is a morphism of Lie-* algebras.
\qed

\sssn
\Cor
For any $D_X$-locally free 
finitely generated Lie-* algebra $\CL$ there exists a central extension
in the class of Lie-* algebras as follows
$$
0\map\Omega_X\map \CL^\flat
\map
\CL\map 0. 
$$
This is just the inverse image of the Tate central extension of $\Mat$.
\qed

Below we denote the Lie-* algebra $\CL^\flat$ by $\CL\Tate$ and cal it {\em
the Tate central extension} of $\CL$.
\vskip 1mm
\noindent
\Rem
Note that for the complete curve $X$ the short exact sequence of the 
$D_X$-modules that defines the extension $\CL\Tate$ does not split, even if we 
forget about the Lie-* algebra structures.

\subsection{Local analog of the Tate extension.}
Fix a point $x\in X$. Let $\hat\CO_x$ (resp. $\hat\CK_x$) be the completion 
of the local ring of the point $x$ (resp. of the local field 
of the point $x$).

Consider the topological Lie algebra  $DR^0(\Spec(\hat\CK_x),\CL)$ with the 
Lie subalgebra $DR^0(\Spec(\hat\CO_x),\CL)$.

\sssn
\Lemma
\begin{itemize}
\item[(i)]
$DR^0(\Spec(\hat\CK_x),\CL)$   is a {\em Tate Lie algebra}.
\item[(ii)]
The subalgebra $DR^0(\Spec(\hat\CO_x),\CL)$ is a {\em c-lattice} in
$DR^0(\Spec(\hat\CK_x),\CL)$.
\qed
\end{itemize}
In particular we have a one dimensional Lie algebra central extension
$$
0\map{\mathbb C}\map DR^0(\Spec(\hat\CK_x),\CL\Tate)\map
DR^0(\Spec(\hat\CK_x),\CL)\map 0.
$$
\sssn
\Prop
The central extension $DR^0(\Spec(\hat\CK_x),\CL\Tate)$ coincides with the 
extension of the Tate Lie algebra $DR^0(\Spec(\hat\CK_x),\CL)$ with the help of 
the {\em critical cocycle} (see e.g. \cite{Ar2}, 4.3.3).
\qed

\section{DG Lie algebroids in the category of $D_X$-modules.}
Recall that $\mod\operatorname{-}D_X$ is a symmetric tensor category.
Below we use this structure to mimick the ordinary definition of a Lie
algebroid in the category of vector spaces.

\ssn
\Def
Let $\CR$ be a $\tenf$-commutative algebra in
$\mod\operatorname{-}D_X$, and let $\CA$ be a $\CR$-module (so we have
a D-module map $\CR\tenf\CA\map\CA$ providing the structure) carrying
a Lie algebra structure in $\mod\operatorname{-}D_X$ (i.e.  a Lie
bracket map $[\cdot,\cdot]:\ \Lambda^2(\CA)\map\CA$ satisfying the
Jacobi identity is given).  Moreover suppose that $\CA$ acts on $\CR$
by derivations (i.e.  we have a D-module map $\CA\tenf\CR\map\CR$ such
that $[a,rb]=a(r)b+r[a,b]$ is satisfied for any sections $a,b\in
\Gamma(X,\CA)$ and $r\in \Gamma(X,\CR)$).

We call the above data {\em the Lie algebroid} in the category 
$\mod\operatorname{-}D_X$ over the $\tenf$-commutative algebra $\CR$.
 
From now
on we assume that all the appearing $\CR$-algebroids are locally free as
$\CR$-modules.

\sssn
By definition a {\em right module}  (resp. a {\em left module}) 
over a Lie algebroid $(\CA,\CR)$ on $X$
is  a sheaf of $\CR$-modules $\CM$ with the Lie
action of $\CA$ satisfying the constraint 
$(rm)\cdot a=r(m\cdot a)-(a(r))m$
(resp.
$a\cdot(rm)=r(a\cdot
m)+(a(r))\cdot m$) 
for any sections $a$ of $\CA$, $r$ of $\CR$ and $m$ of $\CM$.

\sssn
Recall that the {\em universal enveloping algebra} for a $\CR$-Lie
algebroid $\CL$ in $\mod\operatorname{-}D_X$ is defined in the same
way it is done for Lie algebroids over vector spaces:  we take the
free algebra in the category $\mod\operatorname{-}D_X$ generated by
$\CL$ and take its quotient by the obvious ideal of relations
including the one expressing the action of $\CL$ on $\CR$ by
derivations.  We denote the obtained associative algebra in
$\mod\operatorname{-}D_X$ by $\CU_{\CR}(\CL)$.

\subsection{Homological Chevalley complex for a right module over a 
Lie algebroid in \protect{$\mod\operatorname{-}D_X$}.} \label{homol}
For a {\em right} $\CA$-module $\CM$
consider the graded $D_X$-module on $C\bul(\CA,\CM)$ as follows:
$$C\bul(\CA,\CM)=\underset{k}\bigoplus C^k(\CA,\CM),\
C^k(\CA,\CM)=\CM\ten_{\CR}(\Lambda_{\CR}^{-k}(\CA)).$$
We endow the graded vector space with the differential as follows:
For sections $a_1,\ldots,a_p$ of $\CA$ and $m$ of $\CM$ we put
\label{algebroid}
\begin{gather*}
d(m\ten a_1\wedge\ldots\wedge a_p)
=\sum_i(-1)^i
m\cdot a_i\ten a_1\wedge\ldots\wedge \hat a_i\wedge\ldots\wedge a_p\\
+\sum_{i<j}(-1)^{i+j}m\ten [a_i,a_j]\wedge a_1\wedge\ldots\wedge 
\hat a_i\wedge
\ldots\wedge \hat a_j\wedge\ldots\wedge a_p.
\end{gather*}
\label{algebrcor}
\Lemma
\begin{itemize}
\item[(i)] The differential in the complex  is well defined.
\item[(ii)] The differential satisfies $d^2=0$.
\end{itemize}
\dok
(i) Let us perform a calculation showing that the differential
$d:\ C^{-k}\map C^{-k+1}$ is
well defined for $k=2$, the general case is quite similar.
We have
\begin{gather*}
d(m\ten ra_1\wedge a_2)=
m\cdot(ra_1)\ten a_2-
m\cdot a_2\ten ra_1+m\ten [a_2,ra_1]\\=
m\ten (a_2(r))a_1+m\ten r[a_2,a_1]-r(m\cdot a_2)\ten a_1+m\cdot(ra_1)
\ten a_2\\
=(a_2(r))m\ten a_1+rm\ten[a_2,a_1]-(rm)\cdot a_2\ten a_1-
a_2(r)m\ten a_1+(rm)\cdot a_1\ten a_2\\ =
(rm)\cdot a_1\ten a_2
-
(rm)\cdot a_2\ten a_1
+
rm\ten[a_2,a_1]=
d(rm\ten a_1\wedge a_2).
\end{gather*}
The general case is quite similar.

(ii) This is the usual calculation in the Chevalley complex.
\qed
\vskip 1mm
\noindent
\Rem
We have shown above that there exists a {\em homological} Chevalley complex
for a {\em right} $\CA$-module $\CM$. A very similar calculation 
proves the existence
of the {\em cohomological} Chevalley complex for a {\em left} 
$\CA$-module $\CN$, of the size
$\underline{\hom}_{\CR}(\Lambda_{\CR}\bul(\CA),\CN)$.

\sssn 
In fact both the homological and the cohomological versions of the
Chevalley complex for a Lie algebroid appear naturally ``in a
coordinate-free way'' as a result of the following construction.

Consider the tautological {\em left} $\CA$-module $\underline{\CR}$.
We construct its standard projective resolution in the way it is
usually done for modules over Lie algebras.

Namely consider the complex $$
\mathsf{Stan}\bul(\CA, \underline{\CR}):=
\CU_{\CR}(\CA)\ten_{\CR}\Lambda\bul_{\CR}(\CA)\ten_{\CR}\underline{\CR}
$$
with the standard differential. Note that the differential uses the right
$CA$-module structure on $\CU_{\CR}(\CA)$.
\vskip 1mm
\noindent
\Lemma
\begin{itemize}
\item[(i)]
The complex of $D_X$-modules $\mathsf{Stan}\bul(\CA, \underline{\CR})$
is well-defined.
\item[(ii)]
$\mathsf{Stan}\bul(\CA, \underline{\CR})$ is a complex of $\CA$-modules.
\item[(iii)]
The  homological Chevalley complex $C\bul(\CA,\CM)$ for a right $\CA$-module
$\CM$ is isomorphic to 
$\CM\ten_{\CA}\mathsf{Stan}\bul(\CA, \underline{\CR}).
$
\item[(iv)]
The  cohomological Chevalley complex $\underline{\hom}_{\CR}(\Lambda_{\CR}\bul
(\CA),\CN)$
 for a left $\CA$-module
$\CN$ is isomorphic to
$\underline{\hom}_{\CA}(\mathsf{Stan}\bul(\CA, \underline{\CR}),\CN).\qed
$
\end{itemize}
\subsection{DG Lie algebroids in $D_X$-modules.}
Now let $\CR\bul=\oplus_k\CR^k$ be a graded supercommutative
$\tenf$-algebra on $X$ and let $\CA\bul=\oplus \CA^k$ be a graded
super Lie algebroid over $\CR$.  Suppose also there are an odd
derivation $d_{\CR}$ on $\CR\bul$ of degree $1$ and an odd derivation
$d_{\CA}$ of the Lie superalgebra $\CA\bul$ also of degree $1$
satisfying Leibnitz rule with respect one to another.  Moreover both
of them satisfy $d^2=0$.

\sssn
\Def
The data $(\CA\bul,\CR\bul,d_{\CA},d_{\CR})$ are called the 
differential graded Lie algebroid in the category $\mod\operatorname{-}D_X$
 or, for short, a DG Lie algebroid on $X$. 

The notion of a left (resp. right) DG-module over a DG algebroid in 
$\mod\operatorname{-}D_X$
is a natural  combination of the previous definitions and we do not
spell it out explicitly.
The category of left (resp. right) DG-modules over a DG Lie algebroid
$\CA=\CA_X$ is denoted by
$DG\operatorname{-}\CA\bul\operatorname{-}\mod$
(resp.
$DG\operatorname{-}\mod\operatorname{-}\CA\bul$).

\subsection{Homological Chevalley complex for a DG Lie algebroid in 
\protect{$\mod\operatorname{-}D_X$}.}\label{maincomp}
Now we sort of add a second differential on  the homological Chevalley
complex
given in \ref{algebroid}.
For
$\CM\bul\in
DG\operatorname{-}\mod\operatorname{-}\CA\bul$
consider the bigraded vector space
$C^{\bullet\bullet}(\CA\bul,M\bul)$ as follows:
$
C^{\bullet\bullet}(\CA,M\bul)=
\CM\bul\ten_{\CR\bul}(\Lambda_{\CR\bul}\bul(\CA\bul)),
$
here the first grading comes from the number of wedges in the exterior
product and
$$C^{\bullet k}(\CA\bul,\CM\bul)=\left(\CM\bul\ten_{\CR\bul}
(\Lambda_{\CR\bul}\bul(\CA\bul))\right)^k$$
in the graded tensor product sense.

Consider the two differentials on the bigraded vector space.  The
first one of the grading $(1,0)$ is the usual Chevalley differential
like in \ref{algebroid}.  The second differential of the grading
$(0,1)$ is provided by the differentials on $\CM\bul$ and
$\Lambda_{\CR\bul}\bul(\CA\bul)$.

Consider the total grading on the bigraded space and the 
total differential on it.

\sssn
\Lemma
The differential $d_1+d_2$ is well defined and its square equals zero.
\qed
 
\section{Semiinfinite cohomology via DG Lie algebroids in
\protect{$\mod\operatorname{-}D_X$}.}  In this section we show that
the standard complex for the computation of semiinfinite cohomology of
a chiral module over a Lie-* algebra coincides with the homological
Chevalley complex of the form~\ref{maincomp} for a certain DG Lie
algebroid in the category $\mod\operatorname{-}D_X$ and a certain
right module over it.

\subsection{Construction of the Lie algebroid in
\protect{$\mod\operatorname{-}D_X$}.}  Fix a Lie-* algebra $\CL$ on
$X$.  As before, we suppose that it is locally free and finitely
generated over $D_X$.

The construction will be local and we can assume that 
$X$ is affine.

Consider the completion of the $D_{X\times X}$-module $\Omega\tenb
\CL$ along the diagonal.  We denote this D-module by
$\CL_{\hat\Delta}$.  One can view the D-module $\CL_{\hat\Delta}$ as
the restriction of $\Omega\tenb \CL$ to the ``family of formal discs''
parametrized by the diagonal.

Consider also the D-module $j_*(\CO_U)\ten_{\CO_{X\times X}}\CL_{\hat\Delta}$  
denoted by $\CL_{\til\Delta}$. One can view the D-module
$\CL_{\til\Delta}$ as the restriction of $\Omega\tenb \CL$ to 
the ``family of punctured formal discs''
parametrized by the diagonal.  

\sssn
\Lemma
We have the short exact sequence of the d-modules
$$
0\map\CL_{\hat\Delta}\map\CL_{\til\Delta}\map\Delta_!(\CL)\map 0
\qed
$$
Now we take the D-module direct images
$\CB:=p_1{}_{*}\left(\CL_{\hat\Delta}\right)$ and
$\CG:=p_1{}_{*}\left(\CL_{\til\Delta}\right)$ on $X$ (here $p_1$ is
the projection $(x,y)\in X\times X\longrightarrow x\in X$).  There is
an obvious map $\CB\map\CG$ and from the fact that $\CL$ is a Lie-*
algebra we infer that both $\CB$ and $\CG$ are Lie algebras in the
category of right D-modules on $X$.

Note that the stalk of the D-module $\CB$ (resp.  of $\CG$) at a point
$x\in X$ equals $DR^0(\hat\CO_x,\CL)$ (resp.  $DR^0(\hat\CK_x,\CL)$),
where $\hat\CO_x$ (resp.  $\hat\CK_x$) denotes the spectrum of the
completed local ring (resp.  the completed local field) at $x$.  We
abuse some notation here.  \vskip 1mm \noindent \Lemma There exists a
short exact sequence of D-modules on $X$ as follows:  $$ 0\map \CB\map
\CG\map\CL\map 0.  $$ \dok Consider the short exact sequence of
D-modules from the previous Lemma Now take the (D-module) direct image
of the exact sequence under $p_1$.  It is left to the reader that the
sequence obtained as a result coincides with the one we need.  \qed

We choose the basic $\tenf$-supercommutative algebra for our algebroid to be
$\CR\bul(\CL)=\Lambda\bul(\CL^*)$. 

\sssn
\Lemma
\begin{itemize}
\item[(i)]
The Lie algebra in D-modules $\CB$ acts on $\CL$.
\item[(ii)]
The Lie algebra in D-modules $\CB$ acts on the 
$\tenf$-supercommutative algebra $\CR\bul(\CL)$ by derivations.
\end{itemize}
\dok
Note that the second assertion of the Lemma follows from the first 
one since we can extend the action from the generators of $\CR\bul(\CL)$
to the whole algebra by the Leibnitz rule.

Now we construct the action map for (i) explicitly.  We have the
following sequence of morphisms of D-modules
\begin{gather*} 
p_1{}_{*}\left(\CL_{\hat\Delta}\right)\tenf\CL^*\map
p_1{}_{*}p_1^*\left(p_1{}_{*}\left(\CL_{\hat\Delta}\right)\tenf\CL^*\right)
\map p_1{}_{*}\left(\left(\CL_{\hat\Delta}\right)\tenf
p_1^*(\CL^*)\right)\\ \til\map
p_1{}_{*}((\CL^*\tenb\CL)_{\hat\Delta})\overset{\mu}\map
p_1{}_{*}(\Delta_!(\CL^*))\til\map\CL^*.  
\end{gather*} 
Here
$(\CL^*\tenb\CL)_{\hat\Delta}$ denotes the completion of the D-module
$\CL^*\tenb\CL$ along the diagonal in $X\times X$.

Note that the D-module $\Delta_!(\CL^*)$ is locally finite over the
ideal of the diagonal.  Thus the completion of the map
$\CL^*\tenb\CL\map\Delta_!(\CL^*)$ is well defined.

It is left to the reader that the composition of the above maps provides
the action of the Lie algebra $\CB$ on $\CL^*$.
\qed
\vskip 1mm
\noindent
\Cor
The graded D-module $\CA\bul(\CL):=\CB\tenf\CR\bul(\CL)$ carries a natural
structure of a Lie algebroid over the $\tenf$-supercommutative algebra
$\CR\bul(\CL)$.
\qed

\subsection{Construction of the differential on $\CA\bul(\CL)$.}
Note that the differential $d_{\CR}$ on the $\tenf$-supercommutative algebra
$\CR\bul(\CL)$ is already constructed. Moreover it remains to construct 
the component of the differential on $\CB\tenf\CR\bul(\CL)=
\CB\tenf\Lambda\bul(\CL^*)$ as follows: $d_{\CB}:\ \CB\map\CB\tenf\CL^*$.
After that the differential on the whole Lie algebroid is obtained from
$d_{\CR}+d_{\CB}$ by the Leibnitz rule.

Now by Lemma~\ref{hom}(i) we rewrite the map in question as
$\CB\tenb\CL\map\Delta_!(\CB)$ or
$$\left(p_1{}_{*}\left(\CL_{\hat\Delta}\right)\right)\tenb\CL\map
\Delta_!\left(p_1{}_{*}\left(\CL_{\hat\Delta}\right)\right).$$ Note
that $DR^0(\CL)$ acts naturally on every stalk of
$p_1{}_{*}\left(\CL_{\hat\Delta}\right)$ at any point $x\in X$.
Recall that the stalk equals $DR^0(\hat\CO_x,\CL)$.  Thus the sheaf of
Lie algebras acts on it acts on
$\left(p_1{}_{*}\left(\CL_{\hat\Delta}\right)\right)$.  \sssn \Prop
The above action is given by differential operators, i.e.  it lifts to
the required morphism $\CB\tenb\CL\map\Delta_!(\CB)$.

\dok
This is a local assertion. It is  left to the reader to check it.
\qed
\vskip 1mm
\noindent
\Cor
The D-module $\CA\bul(\CL):=\CB\tenf\CR\bul(\CL)$ carries a natural
structure of a DG Lie algebroid over the $\tenf$-supercommutative DG algebra
$\CR\bul(\CL)$.

\sssn
\Rem
In particular the D-module $\CA\bul(\CL\Tate):=
p_1{}_{*}((\CL_{\Tate})_{\hat\Delta})\tenf\Lambda\bul(\CL^*)$ is a DG
Lie algebroid over $\CR\bul(\CL)$.  This follows from the obvious fact
that
$\underline{\hom}_{D_X}(\CL,D_X\ten\Omega_X)=\underline{\hom}_{D_X}
(\CL\Tate,D_X\ten\Omega_X)$.

\subsection{Construction of the left $\CA\bul(\CL)$ DG-module for a
chiral $\CL$-module.}  Note that the Lie algebra in the category of
D-modules $\CG=p_1{}_{*}\left(\CL_{\til\Delta}\right)$ acts naturally
on an arbitrary chiral $\CL$-module $\CM$ as follows:  \begin{gather*}
p_1{}_{*}\left(\CL_{\til\Delta}\right)\tenf\CM\map
p_1{}_{*}\left(\left(\CL_{\til\Delta}\right)\tenf p_1^*(\CM)\right)
\map
p_1{}_{*}\left(\left(\CL_{\til\Delta}\right)\tenf(\CM\tenb\Omega)\right)\\
\map p_1{}_{*}\left((\CM\boxtimes\CL)_{\til\Delta}\right)\overset{\mu}
\map p_1{}_{*}\Delta_!(\CM) \til\map\CM.  \end{gather*} Here
$(\CM\boxtimes\CL)_{\til\Delta}$ denotes the completion of the
D-module $j_*j^*(\CM\boxtimes\CL)$ along the diagonal in $X\times X$.

Note that the D-module
$\Delta_!(\CM)$ is locally finite over the ideal of the diagonal. Thus the 
completion of the 
map $j_*j^*(\CM\tenb\CL)\map\Delta_!(\CM)$  is well defined.

\sssn
For a chiral $\CL$-module $\CM$ consider it as a Lie-* module and 
recall its naive de Rham complex $C\bul(\CL,\CM)=\Lambda\bul(\CL^*)\tenf\CM$.
\vskip 1mm
\noindent 
\Lemma
The complex $C\bul(\CL,\CM)$ has a natural structure of a {\em left} DG-module
over the DG Lie algebroid $\CA\bul(\CL)$.

\dok
The statement of the Lemma follows from the existence of the
$\CG$-module structure on $\CM$ introduced in the beginning of the
present subsection.  \qed

Here we come to the crucial point explaining the phenomenon of the 
Tate extension
in the semiinfinite cohomology of  Lie-* algebras. What we would 
like to do is to consider
the homological Chevalley complex of the DG Lie algebroid 
$(\CA\bul(\CL),\CR\bul(\CL),\ldots)$ with coefficients
in $C\bul(\CL,\CM)$. Yet 
there is no naive way to do it. Somehow we have to make
$M\ten\Lambda\bul(\CL^*)$ into a {\em right} DG-module over our 
DG Lie algebroid in the category of $D_X$-modules. 
\vskip 1mm
\noindent
\Rem
Mimicking the Tate Lie
algebra case we could consider the DG Lie algebroid 
$\CA\bul(\CL)\oplus\CR\bul(\CL)\ten\mathbf{1}$, then construct its antipode 
$\alpha$ 
not commuting with the differential. However the point is that the obtained
DG Lie algebroid in the category $\mod\operatorname{-}D_X$ {\em does not come
from a central extension of our Lie-* algebra} $\CL$. Thus there is no way to 
construct a left DG-module over 
$(\CA\bul(\CL)\oplus\CR\bul(\CL)\ten\mathbf{1})^{\operatorname{opp}}$
(with the differential twisted by the antipode) starting from 
a chiral module either 
over $\CL$ or over some its central extension. 
Instead one should act as follows.

\subsection{Antipode map for the Tate extension of \protect{$\CA\bul(\CL)$}.}
Consider the (Tate) extension of  Lie algebras in D-modules on $X$:
$$
0\map\Omega_X\map\CB_{\Tate}\map\CB\map 0.
$$
Here $\CB_{\Tate}$ denotes $\CB(\CL_{\Tate})$. We have also the 
corresponding extension of DG Lie algebroids on $X$:
$$0\map\Omega_X\tenf\CR\bul(\CL)
\map\CB_{\Tate}\tenf\CR\bul(\CL)\map\CB\tenf\CR\bul(\CL)\map 0.$$
Denote $\CB_{\Tate}\tenf\CR\bul(\CL)$ by $\CA\bul(\CL_{\Tate})$.

We will need also the universal enveloping algebras of the DG-Lie
$\CR\bul(\CL)$-algebroids $\CA\bul(\CL_{\Tate})$ and
$\CA\bul(\CL_{\Tate})$.  Keeping the notation from the previous
section we denote these associative algebras in
$\operatorname{mod-}D_X$ by $\til\CU_{\CR}(\CA\bul(\CL_{\Tate}))$ and
$\CU_{\CR}(\CA\bul(\CL))$ respectively.  Let
$\CU_{\CR}(\CA\bul(\CL_{\Tate}))$ be the qotient of
$\til\CU_{\CR}(\CA\bul(\CL_{\Tate}))$ by the ideal generated by the
relation $\CR\bul(\CL)\tenf(\Omega_X\ten \mathbf{1})=\CR\bul(\CL)$.
Here the LHS of the equality is the kernel of the DG Lie algebroid
extension map while the RHS is the unit
$\CR\bul(\CL)\subset\til\CU_{\CR}(\CA\bul(\CL_{\Tate}))$.

Note that with the differentials forgotten the algebras
$\CU_{\CR}(\CA\bul(\CL_{\Tate}))$ and $\CU_{\CR}(\CA\bul(\CL))$ are
isomorphic.

\sssn
We introduce an antipode map
$\CA\bul(\CL_{\Tate})\til\map\CA\bul(\CL_{\Tate})^{\operatorname{opp}}$
as follows.  Set $\alpha(b\ten r)=-b\ten r
+\operatorname{co-ad}_b(r)\ten \mathbf{1}$.  Here $b$ is a section of
$\CB_{\Tate}$, $r$ is a section of $\CR\bul(\CL)$ and $\mathbf{1}$ is
the generating section of $\Omega_X$.

Note that the antipode does not necesserily commute with the differential on 
$\CA\bul(\CL_{\Tate})$.
\vskip 1mm
\noindent
\Lemma
\begin{itemize}
\item[(i)] $\alpha$ is well defined as a antipode of a Lie superalgebra in the
category of $D_X$-modules $\CA\bul(\CL_{\Tate})$ (with 
the differential forgotten).
\item[(ii)] 
When restricted to any open affine subset $\overset{\circ}{X}\subset X$
the DG Lie superalgebra in the
category of $D_X$-modules $\CA\bul(\CL_{\Tate})^{\operatorname{opp}}
|_{\overset{\circ}{X}}$ with the differential $\alpha\circ d_{\CA\bul}
\circ\alpha^{-1}$
is isomorphic to the DG Lie superalgebra in the
category of $D_X$-modules $\CA\bul(\CL\oplus\Omega_X)|_{\overset{\circ}{X}}$.
Here $\CL\oplus\Omega_X$ denotes the trivial central extension.
\end{itemize}
\dok
Both statements of the Lemma follow from the corresponding local statements 
presented in \cite{Ar2}, Proposition 4.3.4.
\qed

Here we come to another difference with the Tate Lie algebra
semiinfinite cohomology case.  While previous Lemma states that when
restricted to an open affine subset the complex of D-modules
$\CA\bul(\CL_{\Tate})^{\operatorname{opp}}$ is isomorphic to
$\CA\bul(\CL)\oplus\CR\bul(\CL)\tenf (\Omega_X\ten\mathbf{1})$, still
possibly the short exact sequence of D-modules on the complete curve
$$0\map\Omega_X\tenf\CR\bul(\CL)
\map(\CB_{\Tate}\tenf\CR\bul(\CL))^{\operatorname{opp}}
\map\CB\tenf\CR\bul(\CL)\map 0$$ does not split.

That is where we use the universal enveloping algebras of our DG Lie
algebroids.  Extend the antipode $\alpha$ to
$\til\CU_{\CR}(\CA\bul(\CL_{\Tate}))$.

\sssn
\Prop
\begin{itemize}
\item[(i)]
$\alpha$ descends to the antipode of $\CU_{\CR}(\CA\bul(\CL_{\Tate}))$.
\item[(ii)] The DG-algebra 
$\CU_{\CR}(\CA\bul(\CL_{\Tate}))^{\operatorname{opp}}$ with 
the differential twisted by the antipode 
$\alpha$ is isomorphic to the DG-algebra $\CU_{\CR}(\CA(\CL))$.
\end{itemize}

\dok
Follows from local calculations in the Tate Lie algebra semiinfinite cohomology 
case (see \cite{Ar2}).
\vskip 1mm
\noindent
\Cor
Any {\em left} DG-module over the DG Lie algebroid $\CA\bul(\CL_{\Tate})$
on which the center $\Omega_X\ten\mathbf{1}$ acts by 
unity becomes a {\em right}
DG-module over the DG Lie algebroid $\CA\bul(\CL)$.
\qed

\subsection{Standard semiinfinite complex for a chiral module over 
Tate extension of the Lie-* algebra
$\CL$.}
For a chiral $\CL_{\Tate}$-module $\CM$ such that the center $
\Omega_X\ten\mathbf{1}$ acts on it 
by unity consider the  naive de Rham complex of the module
$C\bul(\CL_{\Tate},\CM)=\CM\tenf\Lambda(\CL^*)$ as a {\em right} 
DG-module over the
DG-Lie algebroid in the category of $D_X$-modules 
$(\CB\tenf\CR\bul(\CL))$. The right DG-module structure 
is obtained using the antipode construction from the previous subsection.

\sssn
\Def
We call the homological Chevalley complex of the DG Lie algebroid 
$\left((\CB\tenf\CR\bul(\CL)),
\CR\bul(\CL),\ldots\right)$
with coefficients in the right DG-module $C\bul(\CL_{\Tate},\CM)$
{\em the standard semiinfinite complex} for the chiral module $\CM$
over the Lie-* algebra $\CL_{\Tate}$ and denote it by 
$C^{\si}(\CL_{\Tate},\CM)$.
\vskip 1mm
\noindent
\Rem
Note that as the D-module on $X$ the constructed complex looks as follows:
\begin{gather*}
C^{\si}(\CL_{\Tate},\CM)
=\Lambda\bul_{\CR\bul}\left(p_1{}_{*}((\CL)_{\hat\Delta}
\tenf\CR\bul(\CL)\right)
\ten_{\CR\bul}\left(\CM\tenf\Lambda\bul(\CL^*)\right)\\
=\Lambda\bul\left(p_1{}_{*}((\CL)_{\hat\Delta}\right)\tenf 
\Lambda\bul(\CL^*)\tenf \CM.
\end{gather*}
Here as before $p_1$ denotes the projection $X\times X\map X,\ (x,y)\mapsto x$.

\end{document}